\documentclass{amsart} 

\usepackage{amssymb}
\usepackage{xy} \xyoption{all}

\theoremstyle{plain}
\newtheorem{theorem}{Theorem}[section]

\newtheorem{proposition}[theorem]{Proposition}

\theoremstyle{definition}
\newtheorem{definition}[theorem]{Definition}
\newtheorem{example}[theorem]{Example}

\newtheorem{remark}[theorem]{Remark}

\numberwithin{equation}{section}

\makeatletter

\newcommand{\C}{\mathbb{C}}

\newcommand{\op}{\mathrm{op}}

\DeclareMathOperator{\id}{id}

\newcommand{\BDC}{\mathbf{D}^{\mathrm{b}}}

\newcommand{\DSum}{\bigoplus}

\newcommand{\ilim}[1][]{\mathop{\varinjlim}\limits_{#1}}

\renewcommand{\to}[1][]{\xrightarrow{#1}}
\newcommand{\from}[1][]{\xleftarrow{#1}}
\def\isoto{\@ifnextchar [{\relisoto}{\absisoto}}
\newcommand{\relisoto}[1][]{\xrightarrow[\sim]{#1}}
\newcommand{\absisoto}[1][]{\xrightarrow{\sim}}

\newcommand{\Endo}[1][]{\mathrm{End}_{\raise1.5ex\hbox to.1em{}#1}}

\newcommand{\Hom}[1][]{\mathrm{Hom}_{\raise1.5ex\hbox to.1em{}#1}}

\newcommand{\RHom}[1][]{\mathrm{RHom}_{\raise1.5ex\hbox to.1em{}#1}}

\newcommand{\Ext}[2][]{\mathrm{Ext}_{\raise1.5ex\hbox to.1em{}#1}^{#2}}

\newcommand{\Tens}[1][]{\mathbin{\otimes_{\raise1.5ex\hbox to-.1em{}#1}}}

\newcommand{\LTens}[1][]{\mathbin{\otimes_{\raise1.5ex\hbox to-.1em{}#1}^{L}}}

\newcommand{\Tor}[2][]{\mathrm{Tor}^{\raise1.5ex\hbox to.1em{}#1}_{#2}}

\newcommand{\sheaffont}[1]{\mathcal{#1}}
\def\sha{\sheaffont{A}}

\def\she{\sheaffont{E}}

\def\sho{\sheaffont{O}}

\def\shr{\sheaffont{R}}

\def\shw{\sheaffont{W}}

\renewcommand{\hom}[1][]{{\sheaffont{H}om}_{\raise1.5ex\hbox to.1em{}#1}}
\newcommand{\aut}[1][]{{\sheaffont{A}ut}_{\raise1.5ex\hbox to.1em{}#1}}
\newcommand{\inn}[1][]{{\sheaffont{I}nn}_{\raise1.5ex\hbox to.1em{}#1}}

\newcommand{\rhom}[1][]{{R\sheaffont{H}om}_{\raise1.5ex\hbox to.1em{}#1}}

\newcommand{\ext}[2][]{{\sheaffont{E}xt}_{\raise1.5ex\hbox to.1em{}#1}^{#2}}

\newcommand{\thom}[1][]{{T\sheaffont{H}om}_{\raise1.5ex\hbox to.1em{}#1}}

\newcommand{\tens}[1][]{\mathbin{\otimes_{\raise1.5ex\hbox to-.1em{}#1}}}

\newcommand{\ltens}[1][]{\mathbin{\otimes_{\raise1.5ex\hbox to-.1em{}#1}^{L}}}

\newcommand{\tor}[2][]{{\sheaffont{T}or}^{\raise1.5ex\hbox to.1em{}#1}_{#2}}

\newcommand{\oim}[1]{{#1}_*}

\newcommand{\opb}[1]{#1^{-1}}

\newcommand{\GHom}[1][]{\mathrm{GHom}_{\raise1.5ex\hbox to.1em{}#1}}

\newcommand{\GExt}[2][]{\mathrm{GExt}_{\raise1.5ex\hbox to.1em{}#1}^{#2}}

\newcommand{\FHom}[1][]{\mathrm{FHom}_{\raise1.5ex\hbox to.1em{}#1}}

\newcommand{\ghom}[1][]{{\sheaffont{GH}om}_{\raise1.5ex\hbox to.1em{}#1}}

\newcommand{\gext}[2][]{{\sheaffont{GE}xt}_{\raise1.5ex\hbox to.1em{}#1}^{#2}}

\newcommand{\fhom}[1][]{{\sheaffont{FH}om}_{\raise1.5ex\hbox to.1em{}#1}}

\newcommand{\tenstop}[1][]{\mathbin{\hat{\otimes}_{\raise1.5ex\hbox to-.1em{}#1}}}

\newcommand{\homtop}[1][]{\sheaffont{L}_{\raise1.5ex\hbox to.1em{}#1}}

\newcommand{\Homtop}[1][]{\mathrm{L}_{\raise1.5ex\hbox to.1em{}#1}}

\def\absdoim#1{\underline{#1}_*}
\def\reldoim[#1]#2{\underline{#2}_{|{#1}*}}
\def\doim{\@ifnextchar [{\reldoim}{\absdoim}}

\def\absdeim#1{\underline{#1}_*}
\def\reldeim[#1]#2{\underline{#2}_{|{#1}*}}
\def\deim{\@ifnextchar [{\reldeim}{\absdeim}}

\def\absdopb#1{\underline{#1}^{-1}}
\def\reldopb[#1]#2{\underline{#2}_{|{#1}}^{-1}}
\def\dopb{\@ifnextchar [{\reldopb}{\absdopb}}

\def\absboim#1{\underline{\underline{#1}}_*}
\def\relboim[#1]#2{\underline{\underline{#2}}_{|{#1}*}}
\def\boim{\@ifnextchar [{\relboim}{\absboim}}

\def\absbeim#1{\underline{\underline{#1}}_*}
\def\relbeim[#1]#2{\underline{\underline{#2}}_{|{#1}*}}
\def\beim{\@ifnextchar [{\relbeim}{\absbeim}}

\def\absbopb#1{\underline{\underline{#1}}^*}
\def\relbopb[#1]#2{\underline{\underline{#2}}_{|{#1}}^*}
\def\bopb{\@ifnextchar [{\relbopb}{\absbopb}}

\newcommand{\coh}{\mathrm{coh}}

\newcommand{\reghol}{\mathrm{r-hol}}

\renewcommand{\reghol}{\mathrm{rh}}

\newcommand{\Ga}{\C}

\newcommand{\setdef}{;\ }

\newcommand{\h}{\hbar}

\newcommand{\p}{p}
\newcommand{\q}{q}

\newcommand{\ad}{\operatorname{ad}}
\newcommand{\Ad}{\operatorname{Ad}}

\newcommand{\shHom}[1][]{\sheaffont{H}om_{#1}}

\newcommand{\shEnd}[1][]{\sheaffont{E}nd_{#1}}

\newcommand{\stack}[1]{\mathsf{#1}}

\newcommand{\stka}{\stack{A}}

\newcommand{\stke}{\stack{E}}

\newcommand{\stkr}{\stack{R}}

\newcommand{\stkt}{\stack{T}}
\newcommand{\stkw}{\stack{W}}

\newcommand{\stkquant}{{\widetilde\stke}}
\newcommand{\shquant}{{\widetilde\she}}

\newcommand{\stkMod}[1][]{\stack{Mod}_{#1}}

\newcommand{\stkFun}[1][]{\stack{Fct}_{#1}}

\newcommand{\stktimes}{\circledast}
\newcommand{\stktens}[1][]{\mathbin{\stktimes_{\raise1.5ex\hbox to-.1em{}#1}}}

\makeatother

\title{A note on quantization of complex symplectic manifolds}

\author{Andrea D'Agnolo} 
\author{Masaki Kashiwara} 

\address{Dipartimento di Matematica Pura ed Applicata,
Universit{\`a} di Padova,
via Trieste 63, 35121 Padova, Italy} 

\address{Research Institute for Mathematical Sciences,
Kyoto University,
Kyoto, 606-8502, Japan}



\begin{document}

\maketitle 

\begin{abstract}
To a complex symplectic manifold $X$ we associate a canonical quantization
algebroid $\stkquant_X$. This is modeled on the algebras
$\DSum_{\lambda\in\Ga} \oim\rho \she\, e^{\lambda\h^{-1}}$, where $\rho$ is a
local contactification, $\she$ is an algebra of microdifferential operators
and $\h\in\she$ is such that $\Ad(e^{\lambda\h^{-1}})$ is the
automorphism of $\oim\rho\she$ corresponding to translation by $\lambda$ in
the fibers of $\rho$. Our construction is similar to that of
Polesello-Schapira's deformation-quantization algebroid. The deformation
parameter $\h$ acts on $\stkquant_X$ but is not central. If $X$ is compact,
the bounded derived category of regular holonomic $\stkquant_X$-modules is a
$\C$-linear Calabi-Yau triangulated category of dimension $\dim X+1$.
\end{abstract}

\section*{Introduction}

We construct here a canonical quantization algebroid on a complex symplectic
manifold. Our construction is similar to that of the deformation-quantization
algebroid in~\cite{PS04}, which was in turn based on the construction of the
microdifferential algebroid on a complex contact manifold in~\cite{Kas96}. Let
us briefly recall these constructions.

Let $Y$ be a complex contact manifold. By Darboux theorem, the local model of
$Y$ is an open subset of a projective cotangent bundle $P^*M$. A
microdifferential algebra on an open subset $V\subset Y$ is a $\C$-algebra
locally isomorphic to the ring $\she_M$ of microdifferential operators on
$P^*M$. Let $(\she,*)$ be a microdifferential algebra endowed with an
anti-involution. Any two such pairs $(\she',*')$ and $(\she,*)$ are locally
isomorphic. Such isomorphisms are not unique, and in general it is not
possible to patch the algebras $\she$ together in order to get a globally
defined microdifferential algebra on $Y$. However, the automorphisms of
$(\she,*)$ are all inner and are in bijection with a subgroup of invertible
elements of $\she$. As shown in~\cite{Kas96}, this is enough to prove the
existence of a microdifferential algebroid $\stke_Y$, i.e.\ a $\C$-linear
stack locally represented by microdifferential algebras.

Let $X$ be a complex symplectic manifold. On an open subset $U\subset X$, let
$(\rho,\she,*,\h)$ be a quadruple of a contactification $\rho\colon V\to U$, a
microdifferential algebra $\she$ on $V$, an anti-involution $*$ and an
operator $\h\in\she$ such that $\Ad(e^{\lambda\h^{-1}})$ is the automorphism
of $\oim\rho\she$ corresponding to translation by $\lambda\in\C$ in the fibers
of $\rho$. One could try to mimic the above construction in order to get an
algebroid from the algebras $\oim\rho\she$. This fails because the
automorphisms of $(\rho,\she,*,\h)$ are not all inner, an outer automorphism
being given by $\Ad(e^{\lambda\h^{-1}})$. There are two natural ways out.

The first possibility, utilized in~\cite{PS04}, is to replace the algebra
$\oim\rho\she$ by its subalgebra $\shw$ of operators commuting with $\h$. Then
the action of $\Ad(e^{\lambda\h^{-1}})$ is trivial on $\shw$, and these
algebras patch together to give the deformation-quantization algebroid
$\stkw_X$. This is an alternative construction to that of~\cite{Kon01}, where
the parameter $\h$ is only formal (note however that the methods in loc.\
cit.\ apply to general Poisson manifolds).

The second possibility, which we exploit here, is to make
$\Ad(e^{\lambda\h^{-1}})$ an inner automorphism. This is obtained by replacing
the algebra $\oim\rho\she$ by the algebra $\shquant = \DSum_{\lambda\in\Ga}
\oim\rho \she\, e^{\lambda\h^{-1}}$ (or better, a tempered version of it). We
thus obtain what we call the quantization algebroid $\stkquant_X$, where the
deformation parameter $\h$ is no longer central. The centralizer of $\h$ in
$\stkquant_X$ is equivalent to the twist of $\stkw_X \tens[\C]
(\DSum\nolimits_{\lambda\in\C}\C e^{\lambda\h^{-1}})$ by the gerbe
parameterizing the primitives of the symplectic $2$-form. One should compare
this with the construction in~\cite{GW08}, whose authors advocate the
advantages of quantization (as opposed to deformation-quantization) and of the
complex domain.

There is a natural notion of regular holonomic $\stkquant_X$-module. In fact,
for any Lagrangian subvariety $\Lambda$ of $X$ there is a contactification
$\rho\colon Y\to X$ of a neighborhood of $\Lambda$ in $X$ and a Lagrangian
subvariety $\Gamma$ of $Y$ such that $\rho$ induces a homeomorphism
$\Gamma\isoto\Lambda$. Then, an $\stkquant_X$-module is called regular
holonomic along $\Lambda$ if it is induced by a regular holonomic
$\stke_Y$-module along $\Gamma$.

One of the main features of our construction is that, if $X$ is compact, the
bounded derived category of regular holonomic $\stkquant_X$-modules is a
$\C$-linear Calabi-Yau category of dimension $\dim X+1$.

\section{Stacks and algebroids}\label{se:stacks}

Let us briefly recall the notions of stack and of algebroid (refer
to~\cite{Gir71,Kon01,DP05}).

A prestack $\stka$ on a topological space $X$ is a lax analogue of a presheaf
of categories, in the sense that for a chain of open subsets $W\subset
V\subset U$ the restriction functor $\stka(U)\to \stka(W)$ coincides with the
composition $\stka(U)\to \stka(V)\to \stka(W)$ only up to an invertible
transformation (satisfying a natural cocycle condition for chains of four open
subsets). The prestack $\stka$ is called separated if for any $U\subset X$ and
any $\p,\p'\in\stka(U)$ the presheaf $U\supset V\mapsto
\Hom[\stka(V)](\p|_V,\p'|_V)$ is a sheaf. We denote it by
$\shHom[\stka](\p,\p')$. A stack is a separated prestack satisfying a natural
descent condition.

Let $\shr$ be a commutative sheaf of rings. For $\sha$ an $\shr$-algebra
denote by $\stkMod(\sha)$ the stack of left $\sha$-modules. An $\shr$-linear
stack is a stack $\stka$ such that for any $U\subset X$ and any
$\p,\p'\in\stka(U)$ the sheaves $\shHom[\stka](\p',\p)$ have an
$\shr|_U$-module structure compatible with composition and restriction. The
stack of left $\stka$-modules $\stkMod(\stka) =
\stkFun[\shr](\stka,\stkMod(\shr))$ has $\shr$-linear functors as objects and
transformations of functors as morphisms.

An $\shr$-algebroid $\stka$ is an $\shr$-linear stack which is locally non
empty and locally connected by isomorphisms. Thus, an algebroid is to a sheaf
of algebras what a gerbe is to a sheaf of groups. For $\p\in\stka(U)$ set
$\sha_\p = \shEnd[\stka](\p)$. Then $\stka|_U$ is represented by the
$\shr$-algebra $\sha_\p$, meaning that $\stkMod(\stka|_U) \simeq
\stkMod(\sha_\p)$. An $\shr$-algebroid $\stka$ is called invertible if
$\sha_\p\simeq\shr|_U$ for any $\p\in\stka(U)$.

\section{Quantization of contact manifolds}\label{se:contact}

Let $Y$ be a complex contact manifold. In this section we describe a
construction of the microdifferential algebroid $\stke_Y$ of~\cite{Kas96} and
recall some results on regular holonomic $\stke_Y$-modules.

By Darboux theorem, the local model of
$Y$ is an open subset of the projective cotangent bundle $P^*M$ with $M =
\C^{\frac12(\dim Y+1)}$. By definition, a microdifferential algebra on $Y$ is
a $\C$-algebra locally isomorphic to the ring of microdifferential operators
$\she_M$ on $P^*M$ from~\cite{S-K-K}.

Consider a pair $\p = (\she,\,*)$ of a microdifferential algebra $\she$ on an
open subset $V\subset Y$ and an anti-involution $*$, i.e.\ an isomorphism of
$\C$-algebras $*\colon \she\to\she^\op$ such that $** = \id$. Any two such
pairs $\p'$ and $\p$ are locally isomorphic, meaning that there locally exists
an isomorphism of $\C$-algebras $f\colon\she'\to\she$ such that $f *' = * f$.
Moreover, by~\cite[Lemma~1]{Kas96} the automorphisms of $\p$ are all inner and
locally in bijection with the group
\begin{equation}
\label{eq:bE}
\{b \in \she^\times \setdef b^*b=1,\ \sigma(b)=1 \},
\end{equation}
by $b\mapsto\Ad(b)$. Here $\sigma(b)$ denotes the principal symbol and
$\Ad(b)(a) = aba^{-1}$.

\begin{definition}\label{def:EY}
The microdifferential algebroid $\stke_Y$ is the $\C$-linear stack on $Y$
whose objects on an open subset $V$ are pairs $\p = (\she,\,*)$ as above.
Morphisms $\p'\to\p$ are equivalence classes $[a,f]$ of pairs $(a,f)$ with
$a\in\she$ and $f\colon \p'\isoto \p$. The equivalence relation is given by
$(ab,f) \sim (a,\Ad(b)f)$ for $b$ as in \eqref{eq:bE}. Composition is given by
$[a,f] \circ [a',f'] = [af(a'),ff']$. Linearity is given by $[a_1,f_1] +
[a_2,f_2] = [a_1 + a_2 b,f_1]$ for $b$ as in \eqref{eq:bE} with $f_2f_1^{-1} =
\Ad(b)$.
\end{definition}

\begin{remark}
For $M$ a complex manifold, denote by $\Omega_M$ the sheaf of top-degree forms.
The algebra $\she_{\Omega_M^{1/2}} = \Omega_M^{1/2} \tens[\sho_M] \she_M
\tens[\sho_M] \Omega_M^{-1/2}$ has a canonical anti-involution $*$ given by
the formal adjoint at the level of total symbols. The pair
$(\she_{\Omega_M^{1/2}},*)$ is a global object of $\stke_{P^*M}$ whose sheaf
of endomorphisms is $\she_{\Omega_M^{1/2}}$. Thus 
$\stke_{P^*M}$ is represented by $\she_{\Omega_M^{1/2}}$.
\end{remark}

As the algebroid $\stke_Y$ is locally represented by a microdifferential
algebra, it is natural to consider coherent or regular holonomic
$\stke_Y$-modules. Denote by $\stkMod[\coh](\stke_Y)$ and
$\stkMod[\reghol](\stke_Y)$ the corresponding stacks. For $\Lambda \subset Y$
a Lagrangian subvariety, denote by $\stkMod[\Lambda,\reghol](\stke_Y)$ the
stack of regular holonomic $\stke_X$-modules with support on $\Lambda$. 

Denote by $\C_{\Omega_\Lambda^{1/2}}$ the invertible $\C$-algebroid on
$\Lambda$ such that the twisted sheaf $\Omega_\Lambda^{1/2}$ belongs to
$\stkMod(\C_{\Omega_\Lambda^{1/2}})$.

For an invertible $\C$-algebroid $\stkr$, denote by $\stack{LocSys}(\stkr)$
the full substack of $\stkMod(\stkr)$ whose objects are local systems
(i.e.~have microsupport contained in the zero-section).

By~\cite[Proposition~4]{Kas96} (see also \cite[Corollary~6.4]{DS07}), one has

\begin{proposition}\label{pro:LYsmooth}
For $\Lambda \subset Y$ a smooth Lagrangian submanifold there is an
equivalence
\[
\stkMod[\Lambda,\reghol](\stke_Y) \simeq
\oim{p_1}\stack{LocSys}(\opb{p_1}\C_{\Omega_\Lambda^{1/2}}),
\]
where $p_1\colon \Lambda\times \C^\times\to\Lambda$ is the projection.
\end{proposition}

Recall that a $\C$-linear triangulated category $\stkt$ is called Calabi-Yau
of dimension $d$ if for each $M,N\in\stkt$ the vector spaces
$\Hom[\stkt](M,N)$ are finite dimensional and there are isomorphisms
\[
\Hom[\stkt](M,N)^\vee \simeq \Hom[\stkt](N,M[d]),
\]
where $H^\vee$ denotes the dual of a vector space $H$.

Denote by $\BDC_\reghol(\stke_Y)$ the full triangulated subcategory of the
bounded derived category of $\stke_Y$-modules whose objects have regular
holonomic cohomologies.

The following theorem is obtained in \cite{KS08}\footnote{The statement in \cite[Theorem~9.2~(ii)]{KS08} is not correct. It should be read as Theorem~\ref{th:CY} above} as a corollary of results
from \cite{KK81}.

\begin{theorem}\label{th:CY}
If $Y$ is compact, $\BDC_\reghol(\stke_Y)$ is a $\C$-linear Calabi-Yau
triangulated category of the same dimension as $Y$.
\end{theorem}

\section{Quantization of symplectic manifolds}\label{se:symplectic}

Let $X$ be a complex symplectic manifold. In this section we describe a
construction of the deformation-quantization algebroid $\stkw_X$
of~\cite{PS04}, which we also use to introduce the quantization algebroid
$\stkquant_X$. We then discuss some results on regular holonomic
$\stkquant_X$-modules.

By Darboux theorem, the local model
of $X$ is an open subset of the cotangent bundle $T^*M$ with $M =
\C^{\frac12\dim X}$. A contactification $\rho\colon V\to U$ of an open subset
$U\subset X$ is a principal $\C$-bundle whose local model is the projection
\begin{equation} \label{eq:rho} P^*(M\times\C) \supset \{\tau\neq 0\}\to[\rho]
T^*M \end{equation} given by $\rho(x,t;\xi,\tau) = (x,\xi/\tau)$. Here, the
$\C$-action is given by translation $t\mapsto t+\lambda$. Note that the outer
isomorphism of $\oim\rho\she_{M\times\C}$ given by translation at the level of
total symbols is represented by $\Ad(e^{\lambda\partial_t})$.

Consider a quadruple $\q = (\rho,\she,*,\h)$ of a contactification $\rho\colon
V\to U$, a microdifferential algebra $\she$ on $V$, an anti-involution $*$ and
an operator $\h\in\she$ locally corresponding to $\partial_t^{-1}$. Any two
such quadruples $\q'$ and $\q$ are locally isomorphic, meaning that there
locally exists a pair $\tilde f = (\chi,f)$ of a contact transformation
$\chi\colon \rho' \to \rho$ over $U$ and a $\C$-algebra isomorphism $f\colon
\oim\chi\she'\to \she$ such that $f *' = * f$ and $f(\h') = \h$. Moreover,
by~\cite[Lemma~5.4]{PS04} the automorphisms of $\q$ are locally in bijection
with the group
\begin{equation}
\label{eq:bQ}
\Ga_U \times \{b \in \oim\rho\she^\times\setdef [\h,b]=0,\
b^*b=1,\ \sigma_0(b)=1 \},
\end{equation}
by $(\mu,b)\mapsto(T_\mu,\Ad(be^{\mu\h^{-1}}))$. Here $[\h,b] = \h b - b\h$ is
the commutator and $T_\mu$ denotes the action of $\mu$ on $V$.

Consider the quantization algebra
\[ 
\shquant = \DSum_{\lambda\in\Ga}
\bigl(C_\h^\infty\oim\rho \she\bigr) \, e^{\lambda\h^{-1}}, 
\] 
where $C_\h^\infty\oim\rho \she = \{a \in \oim\rho\she\setdef \ad(\h)^N(a)=0,\
\exists N\geq 0 \}$ locally corresponds to operators in
$\oim\rho\she_{M\times\C}$ whose total symbol is polynomial in $t$. Here
$\ad(\h)(a) = [\h,a]$ and the product in $\shquant$ is given by
\[
(a\cdot e^{\lambda\h^{-1}})(b\cdot e^{\mu\h^{-1}}) =
a\Ad(e^{\lambda\h^{-1}})(b)\cdot e^{(\lambda+\mu)\h^{-1}}.
\]
One checks that $\shquant$ is coherent.

\begin{definition}\label{def:QX}
The quantization algebroid $\stkquant_X$ is the $\C$-linear stack on $X$ whose
objects on an open subset $U$ are quadruples $\q = (\rho,\she,*,\h)$ as above.
Morphisms $\q'\to\q$ are equivalence classes $[\tilde a,\tilde f]$ of pairs
$(\tilde a,\tilde f)$ with $\tilde a\in\shquant$ and $\tilde f\colon \q'\isoto
\q$. The equivalence relation is given by $(\tilde a \tilde b,\tilde f) \sim
(\tilde a,\Ad(\tilde b)\tilde f)$ for $\tilde b = be^{\mu\h^{-1}}$ with
$(\mu,b)$ as in \eqref{eq:bE}. Here $\Ad(\tilde b) = (T_\mu
\chi,\Ad(\tilde b))$. Composition and linearity are given as in
Definition~\ref{def:EY}.
\end{definition}

A similar construction works when replacing the algebra $\shquant$ by its
subalgebra $\shw = C^0_\h\oim\rho\she$ of operators commuting with $\h$.
Locally, this corresponds to operators of $\oim\rho\she_{M\times\C}$ whose
total symbol does not depend on $t$. Then the action of
$\Ad(e^{\lambda\h^{-1}})$ is trivial on $\shw$, and these algebras patch
together to give the deformation-quantization algebroid $\stkw_X$
of~\cite{PS04}.

The parameter $\h$ acts on $\stkquant_X$ but is not central. The centralizer
of $\h$ in $\stkquant_X$ is equivalent to the twist of $\stkw_X \tens[\C]
(\DSum\nolimits_{\lambda\in\C}\C e^{\lambda\h^{-1}})$ by the gerbe
parameterizing the primitives of the symplectic $2$-form.

If $X$ admits a global contactification $\rho\colon Y\to X$ one can construct
as above a $\C$-algebroid $\stke_{[\rho]}$ on $X$ locally represented by
$\oim\rho\she$. Then there are natural functors
\begin{equation}
\label{eq:Erho}
\opb\rho\stke_{[\rho]} \to \stke_Y, \quad
\stke_{[\rho]} \to \stkquant_X.
\end{equation}

\begin{remark}
For $M$ a complex manifold, the algebra $\she_{\Omega_{M\times\C}^{1/2}}$ on
$T^*M$ has an anti-involution $*$ and a section $\h=\partial_t^{-1}$ on the
open subset $\tau\neq 0$. For $\rho$ as in \eqref{eq:rho}, the quadruple
$(\rho,\she_{\Omega_{M\times\C}^{1/2}},*,\h)$ is a global object of
$\stkquant_{T^*M}$ whose sheaf of endomorphisms is
$\shquant_{\Omega_{M\times\C}^{1/2}}$. Thus
$\stkquant_{T^*M}$ is represented by $\shquant_{\Omega_{M\times\C}^{1/2}}$.
\end{remark}

In order to introduce the notion of regular holonomic $\stkquant_X$-modules we
need some geometric preparation.

\begin{proposition}\label{pro:Lambda}
Let $\Lambda$ be a Lagrangian subvariety of $X$. Up to replacing $X$ with an
open neighborhood of $\Lambda$, there exists a unique pair
$(\rho,\Gamma)$ with $\rho\colon Y\to X$ a
contactification and $\Gamma$ a Lagrangian subvariety of $Y$
such that $\rho$ gives a homeomorphism $\Gamma\isoto \Lambda$.
\end{proposition}

Let us give an example that shows how, in general, $\Gamma$ and
$\Lambda$ are not isomorphic as complex spaces.

\begin{example}
Let $X=T^*\C$ with symplectic coordinates $(x;u)$, and let $\Lambda\subset X$
be a parametric curve $\{(x(s),u(s))\setdef s\in\C \}$, with $x(0) = u(0) =
0$. Then $Y = X \times \C$ with extra coordinate $t$, $\rho$
is the first projection and $\Gamma$ is the parametric curve
$\{(x(s),u(s),-f(s))\setdef s\in\C \}$, where $f$ satisfies the equations
$f'(s) = u(s)x'(s)$ and $f(0)=0$. For $x(s) = s^3$, $u(s) = s^7 + s^8$ we have
$f(s) = \tfrac3{10}s^{10} + \tfrac3{11}s^{11}$. This is an example where $f$
cannot be written as an analytic function of $(x,u)$. In fact, $s^{11} =
\frac{11}{10}x(s)u(s) - \frac{110}3 f(s)$ and $s^{11} \notin
\C[\mspace{-1mu}[s^3,s^7+s^8]\mspace{-1mu}]$.
\end{example}

One checks as in~\cite{S-K-K} that the functors induced by \eqref{eq:Erho}
\[
\stkMod[\Gamma,\coh](\stke_{Y}) \from [\Psi]
\stkMod[\Gamma,\coh](\stke_{[\rho]}) \to[\Phi] \stkMod[\Lambda,
\coh](\stkquant_X),
\]
are fully faithful. Denote by $\stkMod[\Gamma,
\reghol](\stke_{[\rho]})$ the full abelian substack of
$\stkMod[\Gamma, \coh](\stke_{[\rho]})$ whose essential image
by $\Phi$ consists of regular holonomic $\stke_{Y}$-modules. Denote by
$\stkMod[\Lambda,\reghol](\stkquant_X)$ the essential image of
$\stkMod[\Gamma, \reghol](\stke_{[\rho]})$ by $\Psi$.

\begin{definition}\label{def:hol}
The stack of regular holonomic $\stkquant_X$-modules is the $\C$-linear
abelian stack defined by
\[
\stkMod[\reghol](\stkquant_X) =
\ilim[\Lambda]\stkMod[\Lambda,\reghol](\stkquant_X).
\]
\end{definition}

As a corollary of Proposition~\ref{pro:LYsmooth} we get

\begin{proposition}\label{pro:Lsmooth}
If $\Lambda \subset X$ is a smooth Lagrangian submanifold, there is an
equivalence
\[
\stkMod[\Lambda,\reghol](\stkquant_X) \simeq
\oim{p_1}\stack{LocSys}(\opb{p_1}\C_{\Omega_\Lambda^{1/2}}),
\]
where $p_1\colon \Lambda\times \C^\times\to\Lambda$ is the projection.
\end{proposition}

\begin{remark}
When $X$ is reduced to a point, the category of regular holonomic
$\stkquant_X$-modules is equivalent to the category of local systems on
$\C^\times$.
\end{remark}

Finally, as a corollary of Theorem~\ref{th:CY} we get

\begin{theorem}
If $X$ is compact, $\BDC_\reghol(\stkquant_X)$ is a $\C$-linear Calabi-Yau
triangulated category of dimension $\dim X+1$.
\end{theorem}

\providecommand{\bysame}{\leavevmode\hbox to3em{\hrulefill}\thinspace}

\end{document}